\documentclass{amsart}
\usepackage{graphicx}
\newtheorem{thm}{Theorem}[section]

\newtheorem{conjecture}[thm]{Conjecture}
\newtheorem{lem}[thm]{Lemma}
\newtheorem{prop}[thm]{Proposition}
\theoremstyle{definition}

\theoremstyle{remark}

\numberwithin{equation}{section}

\begin{document}

\title{Higher Todd Classes and Holomorphic Group Actions }%
\author{Jonathan Block and Shmuel Weinberger }%

\thanks{The authors would like to thanks Jonathan Rosenberg for
enlightening discussions and comments in the preparation of this
article.}

\maketitle {\em to Robert MacPherson with admiration, on the
occasion of his 60th birthday.}

\begin{abstract}  This paper attempts to provide an analogue of the
Novikov conjecture for algebraic (or K\"{a}hler) manifolds.  Inter
alia, we prove a conjecture of Rosenberg's on the birational
invariance of higher Todd genera. We argue that in the algebraic
geometric setting the Novikov philosophy naturally includes
non-birational mappings.
\end{abstract}
\section{Introduction}
This paper describes an attempt to export the Novikov conjecture's
philosophy in the direction of algebraic geometry.  We hope that
our discussion is useful to algebraic geometers and topologists.

Our story begins with Hirzebruch's Riemann-Roch theorem, or even
earlier, with Hirzebruch's signature theorem, which was a lemma in
his proof of Riemann-Roch. Let $M^{4k}$ be a smooth closed
oriented manifold of dimension $4k$. The signature of $M$,
$\mbox{sign}(M)$ is by definition the signature of the symmetric
bilinear form
\[
\cup:  H^{2k}(M ; \mathbb{R})\times H^{2k}(M;\mathbb{R})\to
H^{4k}(M;\mathbb{R})\cong\mathbb{R}
\]
By Poincare duality, this is a nonsingular pairing, and since $2k$
is even, the pairing is symmetric.  Such forms can be
diagonalized, and the signature is the difference in dimensions
between the positive and negative definite parts.

Hirzebruch's signature theorem \cite{Hi} asserts that
\begin{equation}\label{sign_theorem}
\mbox{sign}(M) = \langle L(M), [M]\rangle                   .
\end{equation}
where $L(M)$ is a particular graded polynomial in the Pontrjagin
classes of $M$. Note that as a consequence of this, without
knowing a thing more about $L$, that if $N\to M$ is a finite
cover, then $\mbox{sign}(N) = s\cdot\mbox{sign}(M)$, where $s$ is
the number of sheets of the cover. In fact, Hirzebruch was very
interested in the exact formula for $L$, and it is quite
intimately related to problems as disparate as the number of
exotic differential structures on the sphere \cite{KM} and
Euler-MacLauren formula and lattice point counting problems
\cite{CS1}.

The two sides of the formula \eqref{sign_theorem} are of very
different sorts. The left hand side is, by definition, an oriented
homotopy invariant.  (It is defined cohomologically.)  The right
hand side seems to depend on the smooth structure. Indeed,
ultimately \eqref{sign_theorem} is one of the key ingredients in
Novikov's proof that Pontrjagin classes are (rationally)
topological invariants \cite{No}.  We will return to this later.

It is also quite obvious from Stoke's theorem and
\eqref{sign_theorem} that the right hand side vanishes whenever
$M$ is the boundary of an oriented manifold. (If it merely bounded
a chain, the cohomology class $L(M)$ might not extend). That this
is true for $\mbox{sign}(M)$ was first observed by Thom, and it is
a consequence of Poincare duality for manifolds with boundary.
Hirzebruch's original proof of \eqref{sign_theorem} was a
systematic exploitation of \begin{enumerate}\item  the cobordism
invariance of the signature, \item the multiplicative formula
\[
\mbox{sign}(M\times N) = \mbox{sign}(M)\mbox{sign}(N) \]  and,
\item Thom's calculation of oriented cobordism \cite{Th}.
\end{enumerate}

That $\langle L(M), [M]\rangle$ is a homotopy invariant though,
seems to only follow from the signature theorem; it does not have
an independent explanation.

Are there any other combinations of Pontrjagin classes that
integrate to a homotopy invariant?  This was considered by P.Kahn
in his thesis \cite{Ka}, and there is no other rational homotopy
invariant cobordism invariant.\footnote{Note that all
Stiefel-Whitney numbers are (mod 2) invariants of homotopy type
that are cobordism invariant.  Indeed, it follows from the Wu
formulae  that Stiefel-Whitney classes are themselves homotopy
invariant.}

Novikov, however, suggested that if we give our manifolds
"polarizations" i.e. continuous maps $f: M \to B\pi$, then we can
exploit the fundamental groups to possibly obtain more homotopy
invariants.  (Here $B\pi$ denotes the classifying space of the
group $\pi$; it is an Eilenberg space of type $K(\pi,1)$, a space
with fundamental group $\pi$ and contractible universal
cover.\footnote{Spaces with contractible universal cover are
called aspherical.})

More precisely, Novikov conjectured that if $\alpha$  is an
element of $H^*(B\pi ; \mathbb{Q})$ then $\langle f^*(\alpha) \cup
L(M), [M]\rangle $ is an oriented homotopy invariant.  Subsequent
work of Mischenko and Wall (independently) showed that
polarizations into non-aspherical spaces are useless: only
cohomology classes that come from the induced map on cohomology
from $M \to  B\pi_1(M)$ have any chance of producing (polarized)
homotopy invariants, and as with Kahn's theorem, there is no point
in considering characteristic polynomials other than $L$.

We recommend the survey by Jim Davis \cite{Da} for more about the
foundations of the Novikov conjecture.  (Indeed, there are several
volumes, and many surveys of this problem - we shall avoid the
temptation to give a survey of these surveys, here.)

We shall return to the Novikov conjecture, but it is perhaps not
to soon to expand on the Novikov philosophy in other directions.

The Riemann-Roch theorem of Hirzebruch computes the arithmetic
genus $p_a(M)$ of a projective algerbaic manifold\footnote{Of
course, the index theoretic proof \cite{AS} of the Riemann-Roch
formula removes algebraic hypotheses on M.  We have also ignored,
for now, the extension to $\chi(M, V)$ where $V$ is a holomorphic
vector bundle over M.} in terms of another polynomial, the Todd
class $\mbox{Td}(M)$.
\begin{equation}\label{rr_theorem}
p_a(M) = \langle \mbox{Td}(M), [M]\rangle
\end{equation}

The arithmetic genus is the alternating sum of the dimensions of
the Dolbeault cohomology groups of M,
\begin{equation}\begin{split}
p_a(M)& =\sum_i(-1)^i\mbox{dim}H^{0,i}(M)\\
  &  =
  \sum_i(-1)^i\mbox{dim}H^{i}(M;\mathcal{O}_M)\end{split}\end{equation}
The Todd class is a polynomial in the chern classes of M, and we
shall, again, ignore its precise form.

Now, the analogue of the homotopy invariance consequence of
\eqref{sign_theorem} is that $Td(M)$ integrates to a birational
invariant of a smooth variety $M$.  The reason is that the left
hand side is a birational invariant; indeed, as a consequence of
the Hartog extension theorem, each of the individual cohomology
groups arising in the definition of the arithmetic genus are
birational invariants (see \cite{GrH}).

Then the Novikov philosophy of trying to couple with group
cohomology to extend a general result to a more precise one in the
presence of a fundamental group leads one (and in particular
Rosenberg \cite{Ro}\footnote{Actually, we had made the same
conjecture some years ago, and verified in the case of abelian
fundametantal group using ideas of Lusztig \cite{Lu}, but then
noticed the more general results described below so we let the
matter of finding an analogue of the Novikov conjecture drop.  On
reading \cite{Ro}, we decided to return to the problem.} ) to
conjecture:

\begin{conjecture}(Birational invariance of higher Todd Genera)  If $M$ is a smooth
projective variety and $f:M \to B\pi$ is a continuous map and if
$\alpha$ is an element of $H^*(B\pi ; \mathbb{Q})$ then \[\langle
f^*(\alpha) \cup \mbox{Td}(M), [M]\rangle  \] is a birational
invariant.
\end{conjecture}

Note (see \cite{GrH} or \cite{Ro}) that the (topological)
fundamental group is a birational invariant, so the question makes
sense. Rosenberg in fact showed that for many $\pi$, this is true;
it's true whenever a certain approach to the Novikov conjecture
(the so called "analytic approach") works.  We will see that it's
true in general using resolution of singularities and the
Riemann-Roch theorem of Baum-Fulton-MacPherson
\cite{BFM}\footnote{Rosenberg tells us that he received an email
from Schuermann that this follows from the preprint \cite{BSY}, as
well.}, but again, we are rushing the story. Rosenberg also
observed an analogue of the theorems of Kahn-Mischenko-Wall, that
no other combinations of Chern numbers can be birational
invariant.

This philosophy has had another very notable success, in
differential geometry, regarding the problem of constructing
complete metrics of positive scalar curvature.  In that case, the
"general theorem" is due to Lichnerowicz, as an early consequence
of the Atiyah-Singer index theorem for the Dirac operator (which
also implies Hirzebruch's theorems) \cite{AS}, and asserts:

\begin{thm}(Lichnerowicz's theorem)  If $M$ is a spin manifold with a
metric of positive scalar curvature, then $\langle  \hat{A}(M),
[M] \rangle  = 0$, where $\hat{A}(M)$ is the $\hat{A}$-genus.
\end{thm}

And again, there is the:

\begin{conjecture}(Gromov-Lawson-Rosenberg) If $M$ is a spin manifold
with positive scalar curvature and and $f:M\to  B\pi$ is a
continuous map and if $\alpha$  is an element of $H^*(B\pi ;
\mathbb{Q})$ then $\langle f^*(\alpha ) \cup \hat{A}(M),
[M]\rangle  = 0.$
\end{conjecture}

The analogy between this problem and the Novikov conjecture was
developed in Rosenberg, \cite{Ro1}.  Moreover, see \cite{St} for
an explanation of the "converse theorems", based on a surgery
theorem of Gromov-Lawson and Schoen-Yau, and spin cobordism
calculations of Stoltz, as well as many more positive results.

Having mentioned the Atiyah-Singer theorem and elliptic operators,
it is now inevitable that we bring in K-theory\footnote{  For the
purposes of the Gromov-Lawson-Rosenberg conjecture we really
should introduce real K-theory; although that is a crucial part of
the story, we must suppress it here because it would take us too
far afield.} .  (The Grothendieck part of the story will come in
the next section.)  In the original papers of Atiyah and Singer,
they associated a ``symbol bundle" to any elliptic operator
\footnote{We do not distinguish between an elliptic operator and
an elliptic complex, nor between a sheaf or a complex of sheaves.}
$[D]$ in $K^*(T^*M)$, where $T^*M$ is the cotangent bundle of $M$.
(For us $K^*$ denotes the topological $K$-group.) However, $T^*M$
has a natural symplectic structure, $\omega$ hence is orientable
for K-theory, and we can thus associate to $D$ an element of the
dual homology theory $K_*(M)$. The index theorem then asserts that
$\mbox{ind}(D) = \mbox{dim} \,\,\mbox{ker}(D)
-\mbox{dim}\,\,\mbox{cok} (D) = p_*[D]$ in $\mathbb{Z}$, where $p$
is the constant map from $M\to \mbox{pt}$ a point. We can also
explain the index theorem for families easily in this framework,
but we shall not. (We might suggest that the reader consult
\cite{At} for an early approach to the K-homology class associated
to an elliptic operator, and \cite{HR} for a recent text.)

In all the above examples there are operators, namely the
``signature", Dolbeault, and Dirac operators, which give ``symbol
classes" in $K_*(M)$.  Now, instead of considering $p: M \to
\mbox{pt}$, we consider $f: M \to B\pi$.  This then gives us
$f_*[D] \in K_*(B\pi)$ and one can conjecture appropriate
vanishing or invariance properties of this invariant.  We call
this the {\em integral Novikov conjecture}.

By using the Chern character, rational K-homology is identified
with ordinary rational homology, and the conjectures discussed
above are the vanishing of this homology class by checking that
its pairing with arbitrary cohomology classes vanishes.

Moreover, this inclusion of torsion is extremely significant.  For
instance, if one uses real K-theory, then $KO_*(\mbox{pt})$ has
2-torsion, and one obtains a more general obstruction to positive
scalar curvature \cite{Ht}, which can be used to show that certain
homotopy spheres do not have positive scalar curvature.

Nice as all of this is, it's off in detail.  The integral Novikov
conjecture, as we just stated it, is wrong for signature and
Dirac.  In both cases, it's the same counterexample.  If one
considers lens spaces of high dimensions with fundamental group
$\mathbb{Z}/p$, $p$ a large prime, it is easy to do calculations
to give the non-homotopy invariance of the K-theoretic signature
class. (In defining lens spaces, one takes the quotient of a
sphere under a free linear action of a cyclic group:  by varying
the linear representation of $\mathbb{Z}/p$, one gets many
examples of nondiffeomorphic, but homotopy equivalent manifolds,
see \cite{Mi}.)  Similarly, all lens spaces have positive scalar
curvature, but the Dirac class is nontrivial.  These classes are
all torsion, though, so this issue does not affect the (rational)
Novikov conjecture.

Now, there is a very sensible way to formulate an integral Novikov
conjecture, even in the presence of torsion, and it boils down to
what we said above for $\pi$ torsion free, but we shall not pursue
it here.  (See e.g. \cite{BC}, \cite{We} for some discussions:
essentially one studies the invariants of proper but perhaps
non-free actions on contractible spaces rather just free actions.)

In any case, the next result is not analogous to what occurs for
the other operators: it is too strong.
\begin{prop}\label{highertodd}  The higher Todd K-class $f_*([Dolbeault]) \in
K_*(B\pi)$ is always a birational invariant.
\end{prop}

In the next section we shall prove the proposition above and in
the final section give what we think are examples of phenomena
which more closely follow the Novikov philosophy.
\section{Novikov conjectures and Novikov theorems}
Let us recall the Grothendieck-Riemann-Roch theorem of Baum,
Fulton and MacPherson, \cite{BFM}. For them, the novelty was to
extend the Riemann Roch theorem to singular spaces. We however
only need it for smooth ones, but we like their statement because
it takes values in topological $K$-theory.

Let $K^a_0(X)$ denote the Grothendieck group of coherent sheaves
on the algebraic variety $X$. Grothendieck realized that taking
the Euler characteristic of a coherent sheaf was a special case of
pushing forward in algebraic $K$-theory. Thus the Riemann-Roch
problem is about passing from the algebraic/geometric group
$K^a_0(X)$ to receptacle theories which are presumably easier to
compute with and also understanding how this map behaves with
respect to pushforward. For example, the Todd genus is the
correction needed to make the corresponding map chow groups (or
singular homology) commute, \cite{BS}. Baum, Fulton and MacPherson
emphasized that  there are other interesting targets, \cite{BFM}
they used topological $K$-homology.

Recall how the pushforward in algebraic $K$-theory works. Given a
proper map $f:X\to Y$ between algebraic varieties, one defines
\[ f_*:K_0^a(X)\to K_0^a(Y)
\]
for $A$ a coherent sheaf on $X$ and $[A]$ its class in $K_0^a$
\[
f_*[A]= \sum_i(-1)^i [R^i f_*(A)]
\]
where $R^if_*(A)$ denotes the $i$th higher pushforward. This is
the  sheafification of \[ U\mapsto H^i(f^{-1}(U); A)
\]
It is true that $R^i f_*(A)$ is a coherent sheaf and that the
class of $f_*(A)$ only depends on the class of $[A]$.

To be a suitable receiver, the theory needs to have pushforwards
for proper maps, so that one can compare them with the
pushforwards in algebraic $K$-theory. In topological $K$ theory,
the pushforward is defined using duality and Gysin maps,
\cite{BFM}.

We now recall
\begin{thm}(Baum, Fulton, MacPherson)
In the category of quasi-projective schemes over $\mathbb{C}$,
there is a natural transformation
\[
\alpha:K^a_0(X)\to K_0(X)
\]
Furthermore, $\alpha$ is commutes with proper pushforwards.
\end{thm}

Now for a complete variety $X$ the pushforward of
$[\mathcal{O}_X]\in K^a_0$ to a point is the arithmetic genus. On
the other hand, we can take $\alpha([\mathcal{O}_X])\in K_0(X)$
and then pushforward. This is now some topologically computed
number that equals the arithmetic genus, by the Riemann Roch
theorem. The Todd genus shows up when one further, takes the Chern
character map from $K_0(X)\to H_0(X:\mathbb{Q})$.

Some of the pleasing aspects of the $K_0$-valued Riemann-Roch
theorem are
\begin{enumerate}
\item The Riemann-Roch map is quite easy to define and after the
relevant dualities boils down to the forgetful map on the
cohomological $K$-theories:
\[
K^0_a(X)\to K^0(X)
\]
where $K^0_a(X)$ is the Grothendieck group of algebraic vector
bundles on $X$ and $K^0(X)$ is the $K$-group of topological vector
bundles on $X$ and the map just forgets the algebraic structure.
\item There is no correction term since they are both forms of
$K$-theory. \item $K_0(X)$ can capture torsion information.
\end{enumerate}
We note the following example of pushforward as a rather simple
lemma.
\begin{lem}
Let $f:X\to Y$ be a morphism induced from blowing $Y$ up along a
smooth center. That is $f$ is a blow down morphism. Then
\begin{equation}\begin{split}
R^0f_*(\mathcal{O}_X)=\mathcal{O}_Y \\
R^if_*(\mathcal{O}_X)=0 \mbox{ for } i\ne 0 \end{split}
\end{equation}
That is,
\[ f_*([\mathcal{O}_X])=[\mathcal{O}_Y]
\]
in $K^a_0(Y)$.
\end{lem}
\proof This follows merely from the fact that the fibers of a blow
up over a smooth center are either points or $\mathbb{P}^n$'s and
that in both cases $H^i(-,\mathcal{O})=\mathbb{C}$ if $i=0$ and
$0$ otherwise. \endproof

\proof We now prove Proposition \ref{highertodd}. Like Rosenberg,
\cite{Ro} we will use the weak factorization theorem \cite{AMW}.
\begin{thm}
Let $\varphi:X\dashrightarrow Y$ be a birational map between
complete non-singular algebraic varieties over an algebraically
closed field of characteristic $0$. Let $U$ be an open set where
$\varphi$ is an isomorphism. Then $\varphi$ can be factored into a
sequence of birational maps
\begin{equation}\label{weakfactor}
X\stackrel{\varphi_1}{\dashrightarrow}
X_1\stackrel{\varphi_2}{\dashrightarrow} X_2
\stackrel{\varphi_3}{\dashrightarrow}\cdots
\stackrel{\varphi_{k-1}}{\dashrightarrow}
X_{k-1}\stackrel{\varphi_k}{\dashrightarrow} X_k =Y
\end{equation}
where
\begin{enumerate}
\item $\varphi=\varphi_{k}\circ
\varphi_{k-1}\circ\cdots\circ\varphi_1$ \item $\varphi_i$ are
isomorphisms on $U$ \item either $\varphi_i:X_{i-1}\dashrightarrow
X_i$ or $\varphi_i^{-1}:X_i\dashrightarrow X_{i-1}$ is a morphism
of algebraic varieties (in particular, everywhere defined)
obtained by blowing up a smooth irreducible center disjoint from
$U$.
\end{enumerate}

\end{thm}

Given a birational map $\varphi:X\dashrightarrow Y$, it induces an
isomorphism of fundamental groups
\[
\varphi:\pi_1(X)\to \pi_1(Y)
\]
which we denote simply by $\pi_1$.  Choosing a polarization (in
the sense above,  $\rho: Y\to B\pi_1 $ induces a polarization for
$X$. Now we apply the weak factorization to the birational map
$\varphi$ to factor it as in \eqref{weakfactor}. All the spaces
$X_i$ thus inherit factorizations $\rho_i:X_i\to B\pi_1$ making
all the maps to $B\pi_1$ commute. Thus, in order to prove the
proposition it only remains to show that given a commutative
diagram
\begin{equation}\begin{array}{ccc}
 V& \stackrel{\psi}{\to} & W\\
  \rho_V&  \searrow & \downarrow \rho_W\\
   & & B\pi_1
   \end{array}
   \end{equation}
where $\psi$ is a blowdown morphism along a smooth center, that
\[
\rho_V(\alpha_V(\mathcal{O}_V))=\rho_W(\alpha_W(\mathcal{O}_W)).
\]
Consider the commutative diagram
\begin{equation}\begin{array}{ccccc}
K^a_0(V)  & \stackrel{\alpha_V}{\to} &  K_0(V) & \searrow &  \\
\psi_*\downarrow &     & \psi_*\downarrow&   &    K_0(B\pi_1)   \\
K^a_0(W)     & \stackrel{\alpha_W}{\to} &   K_0(W) &\nearrow &  \\
\end{array}
\end{equation}

According to the lemma above, we know that
\begin{equation}\label{a-1}
\psi_*([\mathcal{O}_V)]=[\mathcal{O}_W]\in K_0^a(W)
\end{equation}
So we have
\begin{equation}\label{a-2}
\rho_W(\alpha_W(\mathcal{O}_W))=
\rho_W(\alpha_W(\psi_*(\mathcal{O}_V))) \end{equation} But by the
Riemann-Roch Theorem this is
\begin{equation}\label{a-3}
\rho_W(\psi_*(\alpha_V(\mathcal{O}_V))). \end{equation} Now since
$\rho_W\circ \psi$ is homotopic to $\rho_V$, we have that
\begin{equation}\label{a-4}
\rho_W(\psi_*(\alpha_V(\mathcal{O}_V)))=\rho_V(\alpha_V(\mathcal{O}_V)))
\end{equation} which finishes the proof. \endproof

If we examine the proof given above we now see why the birational
invariance of the higher Todd class is true generally:  it is
because birational equivalence is {\em hereditary}, that is, it is
a condition that is locally checkable on the image. As such, it is
more closely analogous to Novikov's theorem that rational
Pontrjagin classes are topologically invariant than it is to the
Novikov conjecture\footnote{See \cite{FW}, \cite{ChW} for a
discussion of the Novikov conjecture and Novikov philosophy on
noncompact manifolds (as influenced most directly by Roe and
Higson); in particular the first reference explains how to prove
Novikov's theorem as a consequence of a the Novikov conjecture for
the metric manifold $\mathbb{R}^n$.} .

Let us amplify this point.  If one has a homotopy equivalence $h:
M'\to M$, then one does not at all know that $h$ restricts nicely
to $h^{-1}(U)$ for subsets $U\subset M$.  The Novikov conjecture
actually addresses this.  If $f: M \to S^1$ is a map, then being
able to homotop $h$ so that $h$ restricted to $h^{-1}(N)$, for $N
= f^{-1}(1)$ is a homotopy equivalence, then by the Hirzebruch
signature theorem one would have obtained a proof for the
fundamental class of the circle.  (A sort of converse to this
argument can be given via surgery theory.)  In fact, this is
essentially the method used in \cite{FH} in the first proof for
free abelian groups.   However, it is deep, and requires a
homotopy to see any hereditary aspect.

However, homeomorphisms don't present this problem: they are
hereditary homotopy equivalences!  In fact for all open sets, they
are proper homotopy equivalences.  (In fact, they are bounded
homotopy equivalences on all open subsets, when remetrized to be
complete.)  This is what leads to Novikov's theorem.

In fact, Sullivan realized that all that Novikov used was the
hereditary homotopy equivalence property, and as such applies to
CE maps, i.e. maps with (Cech) contractible point inverses.  This
doesn't give much more, though, because Siebenmann \cite{Si}
showed that all such maps are uniform limits of homeomorphisms,
but it does gain punch if one realizes that as rational homology
is all that's ever used, one gets the same conclusion if the map
were $\mathbb{Q}$-CE, i.e. had (Cech) rationally acyclic point
inverse images.

Again the exact same reasoning shows that if $h: X \to Y$ is a
small resolution (see \cite{GM}) then $h_*(L(X)) = L(Y)$ where $L$
here is the Goresky-MacPherson $L$-homology class of a (suitable)
stratified space.  It is true on the characteristic class level,
because an appropriate statement is true on the sheaf level.  This
same line of thought can lead one to the projection formulae in
\cite{CS2}, which if rephrased purely homologically (rather than
geometrically, as something about ``stratified maps") can be held
to include Novikov's theorem.

As a perhaps surprising negative example, positive scalar
curvature should not be thought of local!  Certainly there is no
local connection between the p.s.c. assumption and characteristic
classes:  \cite{KW} show that any function $\phi: M \to
\mathbb{R}$, which is negative somewhere, e.g. in a little ball,
is the scalar curvature of a metric on a compact manifold M.  The
negative scalar curvature set can be tiny even if there's a "big"
cohomological obstruction to positive scalar
curvature\footnote{This is not the case in the noncompact case,
see \cite{Roe}, \cite{BW}.} .

Moreover, in truth, one cannot get any information from incomplete
metrics - and the restriction of a given metric to an open subset
will be incomplete.  Any manifold has an open subset diffeomorphic
to $T^{n-1}\times\mathbb{R}$, which does not have any complete
p.s.c. metric - so there is no hope of a deformation argument.  In
the end, the p.s.c. condition is more like a global hypothesis
than a local one!

\section{Connections to Group Actions}

In this section we review some, but by no means all,   connections
between the Novikov philosophy and group actions. and use this to
suggest a holomorphic problem which does seem tied to the Novikov
philosophy.  We shall also point out how the known "universal
results", which apply to all fundamental groups are essentially
exploitations of the hereditary nature of hypotheses, and so tend
to be correct for connected groups, but unavailable for
disconnected ones.

Our story here starts off with another theorem about vanishing of
the $\hat{A}$-genus of spin manifolds, here in the presence of a
circle action.

\begin{thm}(Atiyah-Hirzebruch \cite{AH})  If M is a spin manifold
admitting a (nontrivial!) smooth circle action, then $\langle
\hat{A}(M), [M]\rangle  = 0$.
\end{thm}

There is actually a slight connection to Lichnerowicz's result
above.  If one had a compact nonabelian group action, rather than
a circle action, then one can produce \cite{LY} an invariant
positive scalar curvature metric (essentially by making the orbits
have very small diameter).  So, in that case one gets the
Atiyah-Hirzebruch vanishing from the Lichnerowicz.  However, as
there are many manifolds (e.g. tori) with circle actions and no
p.s.c. metrics, the results are quite independent.

It is worth noting, moreover, that the above Atiyah-Hirzebruch
theorem fails for the torsion part of the index of the Dirac
operator.  It is tied just to the rational part.  (Conversely, one
does obtains a deep interesting restriction on the manifolds
admitting smooth nonabelian connected Lie group actions from the
\cite{LY} construction combined with Hitchen's refinement, that
does not seem to have a purely differential topological proof.)

One cannot go too far in guessing a nonsimply connected version of
the above, as the torus has a free circle action, but its higher
signature (associated to the fundamental class) is nontrivial. The
way around this is to note that if a circle acts on any space $X$,
the orbit of the base point defines a loop, whose class $\langle
orbit\rangle \in \pi_1(X)$ is actually central in the fundamental
group. (It is part of the induced map on $\pi_1$ by the map
$S^1\times X \to X$ defining the action!) It thus makes sense to
work with the quotient group $\pi_1(M)/\langle orbit\rangle $.

The following theorem was conjectured by Reinhard Schultz,
motivated by the Novikov and Gromov-Lawson-Rosenberg conjectures.

\begin{thm}(Browder-Hsiang \cite{BH})  If $M$ is a spin manifold,
$S^1$ acts nontrivially on $M$, and $f:M \to B\pi$ classifies the
fundamental group of $M$, then for any   $\alpha\in
H^*(B(\pi/\langle orbit\rangle  ; \mathbb{Q})$ one has $\langle
f^*(\alpha ) \cup \hat{A}(M), [M]\rangle  = 0$.
\end{thm}

The basic idea of their proof is this:  essentially\footnote{We
are oversimplifying for convenience of exposition.} they build an
equivariant map from $M \to B\pi/\langle orbit\rangle $ (where
that latter has a trivial action).  If one believes this, then
theorem follows. Without loss of generality we can think of
$B\pi/\langle orbit\rangle $ as a manifold, and a cohomology class
on it as being dual to a submanifold with stably trivial normal
bundle\footnote{Note that we are working rationally.  The
conventional argument assumes that the cohomology class is odd
dimensional, and one then finds a submanifold with trivial normal
bundle using old results of Serre \cite{Se}.  Then one uses tricks
to reduce to this, because for even dimensional classes cup square
obstructs finding a dual class with trivial normal bundle.
However, even in the even case, there is a dual submanifold whose
normal bundle is stably trivial, and whose Euler class exactly
accounts for the cup square.} . The transverse inverse image of
this submanifold has $\hat{A}$-genus exactly equal to the
associated higher $\hat{A}$-genus of $M$, but this submanifold is
both spin and invariant under the circle action and hence has
vanishing $\hat{A}$-genus.

The proof is thus a perfect example of the locality (not in $M$,
but in the target space $B(\pi/\langle orbit\rangle )$. Somewhat
related, and simpler, is the following result (which does have
some torsion information, but which we shall not explain):

\begin{thm}(\cite{We2})  If $ M$ is an oriented manifold \footnote{Rational
homology manifold actually suffices.} which admits a nontrivial
$S^1$ action with nonempty fixed point set $F$, then there is a
natural orientation on $F$ and the higher signatures of $M$ and
$F$ agree in $H_*(B\pi; \mathbb{Q})$.
\end{thm}

Note that if $F$ is nonempty, $\langle orbit\rangle $ is trivial.
The proof is similarly local over the quotient starting with a
general result: The signature of manifold equals the signature of
its fixed point set.

Irrelevant but irresistible (to us) remark\footnote{This was
observed by Cappell and the second author many years ago as part
of some unwritten joint work.} :  This principle is very useful
for calculations.  As a particular amusing example. if one
considers a toric surface whose moment map has image a triangle in
the plane. This manifold has many different circle actions, as a
torus contains many different circles.  Ignoring the three circles
that define the sides of the triangle, all circles have the same
three fixed points = the vertices of the moment polygon. The sign
of the fixed point is determined by whether the line lies within
the angle at the corresponding fixed point a.k.a. vertex.
Considering first a line at a tiny angle with one side and then
applying the observation that signature must be independent of
this, we see that (almost) every line is actually in exactly one
vertex angle. This is Euclid's theorem that the sum of the angles
of a triangle is 180 degrees, but phrased more to his liking; the
sum of the angles of a triangle is a straight angle.

However, the main point of \cite{We2} and more relevant to us is
that there are versions of the above theorem for certain finite
cyclic group actions that are equivalent to the Novikov
conjecture.

\begin{thm} Suppose $G$ is a finite cyclic group which acts
smoothly and twisted-homologically trivially on a smooth manifold
$M$.  Then the following formula for characteristic classes
\[
f_*(L(M)) = (fi)_*(L(F)\cup k(\nu_F )) \in H_*(B \pi_1(M))
\]
holds if the Novikov conjecture holds for $\pi_1(M)$.  Here $k$ is
a certain characteristic class (discussed below) and $v$ denotes
the equivariant normal to $F$. Conversely, if this formula holds
for all such actions of a specific finite cyclic group, and even
only for free actions, one can deduce the Novikov conjecture for
$\pi_1M$.
\end{thm}

We shall describe the characteristic class in terms that makes
it's version for other elliptic operators transparent.  If $D$ is
an elliptic operator, and $g$ is a self map of $M$ which preserves
$D$, then there is a Lefshetz version of $\mbox{ind}(D)$.  One
considers $L(g, D) = tr(g| \mbox{ker}D) -tr(g | \mbox{cok}D)$.

If $g$ is part of a compact group, then Atiyah and Singer give a
characteristic class formula. $L(g, D)$ is the result of
integrating a characteristic class over $F$, the fixed set of $g$.
For us, $k$ is the result of averaging this local class for all of
the generators of the cyclic group generated by $\langle g\rangle
$. (They all have fixed set $F$, of course.)

Now let us be a little more explicit about the homological
triviality condition.  Firstly we assume that the $G$ action lifts
to the universal cover as part of a $G \times \pi$ action.  (This
is like what happens when one has a circle action with nontrivial
fixed point set.)  Then we can consider the action on the homology
of the universal cover or its compactly supported cohomology or
its cohomology, it doesn't much matter.  We assume this action is
trivial.

That this formula should hold under a suitable Novikov conjecture
hypothesis is most easily seen using the ideas of \cite{RW}, which
gives a Lefshetz point of view on the analytic Novikov ideas.  We
should therefore be a bit more precise about a formal aspect of
the latter now,

At the core of this approach is the theory of $C^*$-algebras.  It
turns out that when a an elliptic operator acts on sections of a
bundle whose fibers are finitely generated projective modules over
a $C^*$-algebra $A$, its kernel and cokernel (after perturbation)
can be thought of as projective modules over that algebra.  Thus
one has an index which lies in $K_0(A)$.  For purposes of the
Novikov conjecture the relevant algebra $C^*\pi$ is a completion
of the integral group ring. One can take an operator on $M$, and
pass to the universal cover and keep track of the $\pi$ action.
Equivalently, one is taking coefficients in the tautological
$C^*\pi$-bundle over $M$ and taking its index.

In any case, one is lead to study a natural map $K_*(B\pi) \to
K_*(C^*\pi)$ which takes the higher symbol index to the difference
between a kernel and a cokernel, i.e. a real live analytic index.
If this map is injective, then we say that the analytic Novikov
conjecture holds, and when the latter holds, one can deduce the
ordinary Novikov conjecture and the Gromov-Lawson-Rosenberg
conjecture.

The homological triviality implies that $g$ acts trivially on
suitable modules, so $tr = \mbox{dimension}$.  So $\mbox{ind}(D)$,
which includes the higher signature, is related to $L(h, D)$ for
each $h$ which generates $\langle g\rangle $.  We thus get many
formulae for $\mbox{ind}(D)$, whose average is the one displayed.

The converse result relies on surgery ideas to construct enough
actions to contradict the formula if the Novikov conjecture fails.
One also needs to give a purely surgery theoretic proof of the
formula on the assumption of the original Novikov conjecture
rather than the analytic one, which is why one needs exactly the
class $k$ rather than any of its non-Galois invariant versions.

Note that this  Lefshetz type localization formula is equivalent
to the Novikov conjecture and thus is at least at this time not
something universally true.  It also doesn't hold integrally (in
general).  The point is that it's hypothesis is not at all local.
The homological triviality hypothesis is global.  When one has a
local reason for homological triviality, then one can indeed prove
a suitable vanishing theorem.  For instance if $M$ had an
equivariant map to $B\pi$ (trivial action) so that all point
inverses were acted on homologically trivially, then indeed one
can directly prove that the higher signatures of $M$ vanish. (See
\cite{We3} for related statements.)

However, we can now state and sketch a holomorphic statement which
is also non-local, and which follows from the Novikov conjecture.

\begin{thm}  Suppose that $M$ is a smooth compact K\"{a}hler manifold and
$G$ is a cyclic group acting holomorphically on $M$ with fixed set
$F$ and trivially on the unreduced holomorphic $L^2$ cohomology of
its universal cover. If the analytic Novikov conjecture for
$\pi_1(M)$ holds, then $f_*(Td(M)) = fi_*(Td(F)\cup k'(\nu_F ))
\in H_*(B\pi_1M)$ for a suitable characteristic class $k'$ of the
equivariant normal bundle of $F$.
\end{thm}

The hypothesis again demands a lift of $G \times \pi$ to the
universal cover of $M$.  The relevant $L^2$ cohomology is the
unreduced one (see e.g. \cite{Fa}) By holomorphic part, we mean
the $0,*$ part in the Hodge decomposition. The condition of
homological triviality means that if we take the differential
forms and compress by the idempotent $p=(1 - 1/|G|\sum  g))$ the
spectrum of the Laplacian does not contain $0$.  This implies the
vanishing of the relevant Lefschetz number is equal to the index
of the Dolbeault.

It is worth noting that the topological hypothesis in the previous
theorem suffices, in light of Hodge theory to imply the hypothesis
of this theorem. (See \cite{Gr}, \cite{Ste} for some discussion of
Hodge theory on the universal covers of K\"{a}hler manifolds.)
The reason is this:  If the cohomology of the universal cover
(when contracted by the nontrivial action idempotent) then on the
chain level this complex is acyclic.  Since $L^2(\pi)$ is a module
over $\mathbb{Q}\pi)$, a fortiori that cohomology vanishes as
well. (Again it might be helpful to consult \cite{Fa} or
\cite{Lue}).

Note now that not only is an edge of the Hodge diamond killed by
this argument, but in fact the whole thing is.  The $(k, l)$ part
is also Dolbeault with coefficients in a bundle and its index is
can be computed by Riemann-Roch as well.  Indeed, Hirzebruch's
$\chi_y$-genus encodes all of these simultaneously as a formal
polynomial.  In fact, Hirzebruch used $\chi$ to give a
simultaneous description of Riemann-Roch and the signature formula
(for smooth varieties) in the course of proving the latter.  His
theorem is that
\begin{equation}
\chi_y(M) = \langle T_y(M), [M] \rangle
\end{equation}
for a suitable characteristic class.  Combining all of the above
we obtain:

\begin{thm} Suppose that $M$ is a smooth algebraic variety and $G$
is a cyclic group acting homologrphically on $M$ with fixed set
$F$ and trivially on the cohomology with compact supports of its
universal cover, then if the analytic Novikov conjecture for
$\pi_1M$ holds, then $f_*( T_y (M)) = fi_*( T_y (F)\cup k_y'(\nu_F
)) \in H_*(B\pi_1M)$ for a suitable characteristic class $k_y'$ of
the equivariant normal bundle of $F$ (produced from the
Atiyah-Singer integrand by the usual averaging procedure).
\end{thm}

    We close with a number of problems.
\begin{enumerate}
\item  The results about circle actions held for both the
$L$-genus and the $\hat{A}$-genus.  Presumably there is some
simultaneous generalization that applies to higher elliptic
genera.

\item  Are there examples that show that our vanishing theorem for
higher Todd genera are false for finite $\pi$, if one does not
rationalize.

\item   If $f: M' \to M$ is a holomorphic map which induces an
isomorphism on $H_*$ for $\pi$ covers, do they have the same
higher $T_y$-genus in $H_*(B\pi)$?  For the ordinary Novikov
conjecture, such a generalization is possible, and is part of the
proof of the higher Lefshetz localization theorem.

\item   Is there a version of Nielson theory for elliptic
operators other than De Rham?
\end{enumerate}

\noindent Department of Mathematics University of Pennsylvania
Philadelphia, PA, USA blockj@math.upenn.edu \\

\noindent Department of Mathematics University of Chicago Chicago,
IL, USA \newline shmuel@math.uchicago.edu \newline  currently
visiting Courant Institute of the Mathematical Sciences New York,
NY USA


\begin{thebibliography}{99}

\bibitem{AMW} D. Abramovich, K. Karu, K. Matsuki, J. W{\l}odarczyk, Torification
and factorization of birational maps, {\em J. Amer. Math. Soc.},
vol.  {\bf 15}, 2002, No. 3, 531-572.

\bibitem{At} M. Atiyah, Global theory of elliptic operators. 1970
Proc. Internat. Conf. on Functional Analysis and Related Topics (Tokyo, 1969)
pp. 21--30 Univ. of Tokyo Press, Tokyo

\bibitem{AH} M.Atiyah and F.Hirzebruch, Spin-manifolds and
group actions. 1970  Essays on Topology and Related Topics
(Mémoires dédiés à Georges de Rham)  pp. 18--28 Springer, New York

\bibitem{AS} M.Atiyah and I.Singer, The index of elliptic operators. III.  Ann.
of Math. (2) 87 1968 546--604.

\bibitem{BC} P.Baum and A.Connes,
Geometric $K$-theory for Lie groups and foliations.  Enseign.
Math. (2)  46 (2000), no. 1-2, 3--42.



\bibitem{BFM} P.Baum, W.Fulton, and R.MacPherson, Riemann-Roch and
topological $K$ theory for singular varieties. Acta Math.  143
(1979), no. 3-4, 155--192.

\bibitem{BW} J.Block and S.Weinberger, Large scale homology theories and
geometry. Geometric topology (Athens, GA, 1993),  522--569, AMS/IP
Stud. Adv. Math., 2.1, Amer. Math. Soc., Providence, RI, 1997.

\bibitem{BS} A.Borel and J.P.Serre, Le théorème de Riemann-Roch. (French)
Bull. Soc. Math. France  86 1958 97--136.

\bibitem{BH}  W. Browder, W. C. Hsiang, $G$-actions and the fundamental group,
Invent. Math.  65  (1981/82), no. 3, 411--424.

\bibitem{BSY} P.Brasselet, J. Schuermann, and S.Yokura, Hirzebruch classes and
motivic Chern classes for singular spaces. math.AG/0503492

\bibitem{CS1} S.Cappell and J.Shaneson, Euler-Maclaurin expansions for
lattices above dimension one.  C. R. Acad. Sci. Paris Sér. I Math.
321 (1995), no. 7, 885--890.

\bibitem{CS2} ------- and --------, Stratifiable maps and topological
invariants.  J. Amer. Math. Soc. 4  (1991), no. 3, 521--551.

\bibitem{ChW} S.Chang and S.Weinberger, On Novikov type conjectures, Clay
Institute volume on Noncommutative geometry, ed by N.Higson,
J.Roe, and J.Rosenberg, to appear.

\bibitem{Da} J.Davis, Manifold aspects of the Novikov conjecture.  Surveys
on surgery theory, Vol. 1, 195--224, Ann. of Math. Stud., 145,
Princeton Univ. Press, Princeton, NJ, 2000.

\bibitem{Fa} M.Farber, Homological algebra of Novikov-Shubin invariants
and Morse inequalities. (English. English summary) Geom. Funct.
Anal. 6 (1996), no. 4, 628--665.

\bibitem{FH} F.T.Farrell and W.C.Hsiang, Manifolds with $\pi
\sb{1}=Z\sp{k}$. Manifolds---Amsterdam 1970 (Proc. Nuffic Summer
School),  pp. 36--43. Lecture Notes in Math., Vol. 197, Springer,
Berlin, 1971.

\bibitem{FW} S.Ferry and S.Weinberger, A coarse approach to the Novikov
conjecture.  Novikov conjectures, index theorems and rigidity,
Vol. 1 (Oberwolfach, 1993), 147--163, London Math. Soc. Lecture
Note Ser., 226, Cambridge Univ. Press, Cambridge, 1995..

\bibitem{GM} M.Goresky and R.MacPherson, Intersection homology theory.
Topology 19  (1980), no. 2, 135--162. Part II. Invent. Math.  72
(1983), no. 1, 77--129.

\bibitem{GrH} P.Griffiths and J.Harris, {\em Principles of algebraic geometry},
Pure and Applied Mathematics. Wiley-Interscience, John Wiley and
Sons, New York, 1978

\bibitem{Gr} M.Gromov, Kähler hyperbolicity
and $L\sb 2$-Hodge theory.  J. Differential Geom.  33  (1991),
no. 1, 263--292.

\bibitem{GL} M.Gromov, and B.Lawson,
Positive scalar curvature and the Dirac operator on complete
Riemannian manifolds. Inst. Hautes Études Sci. Publ. Math.  No. 58
(1983), 83--196 (1984)

\bibitem{FRR} S.Ferry, A.Ranicki, and J.Rosenberg, editors, Novikov
conjectures, index theorems and rigidity. Vols. 1 and 2. Including
papers from the conference held at the Mathematisches
Forschungsinstitut Oberwolfach, Oberwolfach, September 6--10,
1993. Edited by Steven C. Ferry, Andrew Ranicki and Jonathan
Rosenberg. London Mathematical Society Lecture Note Series, 226.
Cambridge University Press, Cambridge, 1995.

\bibitem{HR} N.Higson and J.Roe, Analytic $K$-homology. Oxford
Mathematical Monographs. Oxford Science Publications. Oxford
University Press, Oxford, 2000.

\bibitem{Hi} F.Hirzebruch, Topological methods in algebraic geometry.
Third enlarged edition. New appendix and translation from the
second German edition by R. L. E. Schwarzenberger, with an
additional section by A. Borel. Die Grundlehren der Mathematischen
Wissenschaften, Band 131 Springer-Verlag New York, Inc., New York
1966

\bibitem{Ht} N.Hitchin, Harmonic spinors.  Advances in Math.  14 (1974),
1--55.

\bibitem{Ka} P.Kahn, Characteristic numbers and oriented homotopy type.
Topology  3  1965 81--95.

\bibitem{KM} M.Kervaire and J.Milnor, Bernoulli numbers, homotopy groups,
and a theorem of Rohlin. 1960 Proc. Internat. Congress Math. 1958
pp. 454--458 Cambridge Univ. Press, New York

\bibitem{KW} J. Kazdan, F. W. Warner, Existence and conformal
deformation of metrics with prescribed Gaussian and scalar
curvatures,
 Ann. of Math. (2)  101  (1975), 317--331.

\bibitem{LY} B.Lawson and S.T.Yau, Scalar curvature, non-abelian group
actions, and the degree of symmetry of exotic spheres. Comment.
Math. Helv.  49 (1974), 232--244.

\bibitem{Lu} G.Lusztig, Novikov's higher signature and families of
elliptic operators.  J. Differential Geometry  7 (1972), 229--256.


\bibitem{Lue} W.Lueck, $L\sp 2$-invariants: theory and applications to
geometry and $K$-theory. Ergebnisse der Mathematik und ihrer
Grenzgebiete. 3. Folge. A Series of Modern Surveys in Mathematics
[Results in Mathematics and Related Areas. 3rd Series. A Series of
Modern Surveys in Mathematics], 44. Springer-Verlag, Berlin, 2002.


\bibitem{Mi} J. Milnor, Whitehead torsion. Bull. Amer. Math. Soc.  72 1966
358--426.

\bibitem{No} S.Novikov, Algebraic construction and properties of Hermitian
analogs of $K$-theory over rings with involution from the
viewpoint of Hamiltonian formalism. Applications to differential
topology and the theory of characteristic classes. I. II.  Math.
USSR-Izv.  4 (1970), 257--292; ibid. 4 (1970), 479--505;

\bibitem{Ro} J.Rosenberg, An analogue of the Novikov Conjecture in complex
algebraic geometry. math.AG/0509526

\bibitem{Ro1} J. Rosenberg, $C\sp{*} $-algebras, positive scalar
curvature, and the Novikov conjecture. Inst. Hautes Études Sci.
Publ. Math.  No. 58 (1983), 197--212 (1984).

\bibitem{RW} ----- and S.Weinberger, Higher $G$-signatures for Lipschitz
manifolds. $K$-Theory  7 (1993),  no. 2, 101--132.

\bibitem{Roe} J.Roe, An index theorem on open manifolds. I, II.  J.
Differential Geom.  27 (1988),  no. 1, 87--113, 115--136.

\bibitem{Si} L.Siebenmann, Approximating cellular maps by homeomorphisms.
Topology  11 (1972), 271--294

\bibitem{Se} J.P.Serre, Homologie singulière des espaces fibrés.
Applications. (French) Ann. of Math. (2)  54, (1951). 425--505.

\bibitem{Ste} M.Stern, Index theory for certain complete Kähler manifolds.
J. Differential Geom.  37 (1993),  no. 3, 467--503.

\bibitem{St} S.Stoltz, Positive scalar curvature metrics---existence and
classification questions. Proceedings of the International
Congress of Mathematicians, Vol. 1, 2 (Zürich, 1994),  625--636,
Birkhäuser, Basel, 1995.

\bibitem{Th} R.Thom, Quelques propri\'{e}t\'{e}s globales des variétés
différentiables. (French) Comment. Math. Helv.  28, (1954).
17--86.

\bibitem{We} S.Weinberger Aspects of the Novikov conjecture. Geometric and
topological invariants of elliptic operators (Brunswick, ME,
1988),  281--297, Contemp. Math., 105, Amer. Math. Soc.,
Providence, RI, 1990

\bibitem{We2}
S. Weinberger, Group actions and higher signatures.I, Proc. Nat.
Acad. Sci. U.S.A.  82  (1985),  no. 5, 1297--1298. II.  Comm. Pure
Appl. Math.  40  (1987),  no. 2, 179--187.

\bibitem{We3} S. Weinberger, Class
numbers, the Novikov conjecture, and transformation groups.
Topology  27 (1988),  no. 3, 353--365.







\end{thebibliography}
\end{document}